\documentclass[aop,MSNbibl,seceqn,number,citesort,dvips]{arximspdf}

\doi{10.1214/11-AOP719}
\volume{41}
\issue{3B}
\pubyear{2013}
\firstpage{2047}
\lastpage{2065}

\makeatletter
\newcommand{\C}{\mathbf{C}}
\newcommand{\R}{\mathbf{R}}
\newcommand{\pr}{\mathbf{P}}
\newcommand{\ex}{\mathbf{E}}
\newcommand{\ind}{\mathbf{1}}
\newcommand{\eps}{\varepsilon}
\newcommand{\ph}{\varphi}
\newcommand{\ro}{\varrho}
\newcommand{\catalan}{\mathcal{C}}

\newcommand{\laplace}{\mathcal{L}}

\newtheorem{theorem}{Theorem}[section]
\newtheorem{corollary}[theorem]{Corollary}
\newtheorem{proposition}[theorem]{Proposition}
\newproclaim{remark}[theorem]{Remark}
\newproclaim{notation}{Notation}

\newcommand{\re}{\operatorname{Re}}
\newcommand{\im}{\operatorname{Im}}
\newcommand{\var}{\operatorname{Var}}

\newcommand{\eqref}[1]{(\ref{#1})}
\newcommand{\fracf}[2]{({#1})/({#2})}
\newcommand{\fraca}[2]{{#1}/{#2}}
\makeatother

\begin{document}
\begin{frontmatter}

\title{Suprema of L{\'e}vy processes\thanksref{TT1}}
\runtitle{Suprema of L{\'e}vy processes}

\begin{aug}
\author[A]{\fnms{Mateusz} \snm{Kwa{\'s}nicki}\corref{}\ead[label=e1]{mateusz.kwasnicki@pwr.wroc.pl}\thanksref{au1}},
\author[B]{\fnms{Jacek} \snm{Ma{\l}ecki}\ead[label=e2]{jacek.malecki@pwr.wroc.pl}}
\and
\author[C]{\fnms{Micha{\l}} \snm{Ryznar}\ead[label=e3]{michal.ryznar@pwr.wroc.pl}}
\runauthor{M. Kwa{\'s}nicki, J. Ma{\l}ecki and M. Ryznar}
\affiliation{Wroc{\l}aw University of Technology and Polish Academy
of Sciences,
Wroc{\l}aw University of Technology and Universit\'e d'Angers, and
Wroc{\l}aw~University~of~Technology}

\address[A]{M. Kwa{\'s}nicki\\
Institute of Mathematics\\ \quad  and Computer Science\\
Wroc{\l}aw University of Technology\\
ul. Wybrze{\.z}e Wyspia{\'n}skiego 27\\
50-370 Wroc{\l}aw\\ Poland\\
and\\
Institute of Mathematics\\ Polish Academy of Sciences\\ ul. {\'
S}niadeckich 8\\ 00-976 Warszawa\\
Poland\\
\printead{e1}} 

\address[B]{J. Ma{\l}ecki\\
Institute of Mathematics\\ \quad  and Computer Science\\
Wroc{\l}aw University of Technology\\
ul. Wybrze{\.z}e Wyspia{\'n}skiego 27\\
50-370 Wroc{\l}aw\\ Poland\\
and\\
LAREMA \\ Universit\'e d'Angers \\ 2 Bd Lavoisier \\ 49045 Angers cedex
1\\
France\\
\printead{e2}}

\address[C]{M. Ryznar\\
Institute of Mathematics\\ \quad  and Computer Science\\
Wroc{\l}aw University of Technology\\
ul. Wybrze{\.z}e Wyspia{\'n}skiego 27\\
50-370 Wroc{\l}aw\\ Poland\\
\printead{e3}}
\end{aug}
\thankstext{TT1}{Supported by the Polish Ministry of Science and
Higher Education Grant N~N201 373136.}
\thankstext{au1}{Supported by the Foundation for Polish Science.}

\received{\smonth{3} \syear{2011}}
\revised{\smonth{9} \syear{2011}}

%
\begin{abstract}
In this paper we study the supremum functional $M_t = \sup_{0 \le s
\le t} X_s$, where $X_t$, $t \ge0$, is
a one-dimensional L{\'e}vy process. Under very mild assumptions we
provide a simple, uniform estimate of the
cumulative distribution function of $M_t$. In the symmetric case we
find an integral representation of the
Laplace transform of the distribution of $M_t$ if the L{\'
e}vy--Khintchin exponent of the process increases on $(0,\infty)$.
\end{abstract}

%

\end{frontmatter}

\section{Introduction}
\label{secintro}

By a classical reflection argument, the supremum functional $M_t = \sup
_{0 \le s \le t} X_s$ of the Brownian motion $X_t$ has truncated normal
distribution, $\pr(M_t \ge x) = 2 \pr(X_t \ge x)$ ($x \ge0$). A
similar question for symmetric $\alpha$-stable processes was first
studied by Darling~\cite{bibd56}, and the case of general L{\'e}vy
processes $X_t$ was addressed by Baxter and Donsker~\cite{bibbd57}.
Theorem~1 therein gives a formula for the double Laplace transform of
the distribution of $M_t$, which for a \emph{symmetric} L{\'e}vy
process $X_t$ with L{\'e}vy--Khintchin exponent $\Psi(\xi)$ reads
%
\begin{eqnarray}\label{eqbd}
&& \int_0^\infty\int_0^\infty e^{-\xi x - z t} \pr(M_t \in dx)\,dt
\nonumber
\\[-8pt]
\\[-8pt]
&& \qquad = \frac{1}{\sqrt{z}} \exp\biggl( -\frac{1}{\pi} \int
_0^\infty \frac{\xi\log(z + \Psi(\zeta))}{\xi^2 + \zeta^2}\,
d\zeta \biggr) .
\nonumber
\end{eqnarray}
Inversion of the double Laplace transform is typically a very difficult
task. Apart from the Brownian motion case, an explicit formula for the
distribution of $M_t$ was found for the Cauchy process (the symmetric
$1$-stable process) by Darling~\cite{bibd56}, for a compound Poisson
process with $\Psi(\xi) = 1 - \cos\xi$ by Baxter and Donsker \cite
{bibbd57} and for the Poisson process with drift by Pyke~\cite{bibp59}.

The development of the fluctuation theory for L{\'e}vy processes
resulted in many new identities involving the supremum functional
$M_t$; see, for example,~\cite{bibb96,bibd07,bibk06,bibs99}. There are
numerous other representations for the distribution of $M_t$, at least
in the stable case; see \cite
{bibbdp08,bibb73,bibd56,bibd87,bibgj09,bibgj10,bibhk11,bibhk09,bibk10a,bibk10,bibmk10,bibz64}.
The main goal of this article is to give a more explicit formula for
$\pr(M_t < x)$ and simple sharp bounds for $\pr(M_t < x)$ in terms of
the L{\'e}vy--Khintchin exponent $\Psi(\xi)$ for a class of L{\'e}vy
processes. Most estimates of the cumulative distribution function of
$M_t$ are proved for very general L{\'e}vy processes, without symmetry
assumptions.

Let $\tau_x$ denote the first passage time through a barrier at the
level $x$ for the process $X_t$,
\[
\tau_x = \inf\{ t \ge0 \dvtx X_t \ge x \} , \qquad x \ge0 ,
\]
with the infimum understood to be infinity when the set is empty. We
always assume that $X_0 = 0$. Since $\pr(M_t < x) = \pr(\tau_x >
t)$, the problems of finding the cumulative distribution functions of
$M_t$ and $\tau_x$ are the same. The supremum functional and first
passage time statistics are important in various areas of applied
probability~\cite{bibak05,bibbnmr01}, as well as in mathematical
physics~\cite{bibka08,bibkkm07}. The recent progress in the
potential theory of L\'evy processes is, in part, due to the
application of fluctuation theory; see \cite
{bibbbkrsv09,bibcks11,bibgr08,bibksv9,bibksv10a,bibksv10,bibksv11}.

The paper is organized as follows. Section~\ref{secpre} contains some
preliminary material related to Bernstein functions, Stieltjes
functions and estimates for the Laplace transform. In Section~\ref
{secsup} (Theorem~\ref{thsupest} and Corollary~\ref
{corkapparegular1}) we prove, under mild assumptions, the estimate
\[
\pr(M_t < x) \approx\min\bigl( 1, {\kappa(1/t,0) V(x)} \bigr) , \qquad t, x
> 0 ,
\]
where $V(x)$ and $\kappa(z, 0)$ are the renewal function for the
ascending ladder-height process, and the Laplace exponent of the the
ascending ladder-time process corresponding to $X_t$, respectively.
Here $f(x) \approx g(x)$ means that there are constants $c_1, c_2 > 0$
such that $c_1 g(x) \le f(x) \le c_2 g(x)$. In Section~\ref{secv} we
show that in the symmetric case, given some regularity of $\Psi(\xi
)$, we have
\[
V(x) \approx\frac{1}{\sqrt{\Psi(1/x)}} , \qquad x > 0;
\]
see Theorem~\ref{thvest}. Therefore the estimate of the above
cumulative distribution function of $M_t$ takes a very explicit form,
\[
\pr(M_t < x) \approx\min\biggl( 1, \frac{1}{\sqrt{t \Psi(1/x)}} \biggr)
, \qquad t, x > 0 .
\]
The other main result of Section~\ref{secv} is an explicit formula for
the (single, in the space variable) Laplace transform of the
distribution of $M_t$ (Theorem~\ref{thsuplaplace}), under the
assumption that $X_t$ is symmetric and $\Psi(\xi)$ is increasing on
$[0, \infty)$.

When $\Psi(\xi) = \psi(\xi^2)$ for a complete Bernstein function
$\psi(\xi)$, the above results can be significantly improved.
Following the approach of~\cite{bibmk10}, a (rather complicated)
explicit formula for $\pr(M_t < x)$ can be given, and estimates and
asymptotic formulae for $\pr(M_t < x)$ extend to $(d / dt)^n \pr(M_t
< x)$ when $x$ is small or $t$ is large. These results will be covered
in a forthcoming paper.

\begin{notation*}
We denote by $C$, $C_1$, $C_2$, etc. constants in theorems, and by $c$,
$c_1$, $c_2$, etc. temporary constants in proofs. Any dependence of a
constant on some parameters is always indicated by writing, for
example, $c(n, \eps)$. We write $f(x) \sim g(x)$ when $f(x) / g(x) \to
1$. We use the terms \emph{increasing}, \emph{decreasing}, \emph
{concave}, \emph{convex function}, etc. in the weak sense.
\end{notation*}

\section{Preliminaries}
\label{secpre}

\subsection{Complete Bernstein and Stieltjes functions}
\label{seccbf}

A function $\psi(\xi)$ is said to be a \emph{complete Bernstein
function} (CBF) if
%
\begin{equation}\label{eqcbf}
\psi(\xi) = c_1 + c_2 \xi+ \frac{1}{\pi} \int_{0+}^\infty\frac
{\xi}{\xi+ \zeta} \frac{\mu(d\zeta)}{\zeta} , \qquad \xi\in\C
\setminus(-\infty, 0) ,
\end{equation}
where $c_1, c_2 \ge0$, and $\mu$ is a measure on $(0, \infty)$ such
that the integral\break $\int_0^\infty\min(\zeta^{-1},  \zeta^{-2}) \mu
(d\zeta)$ is finite. A function $\tilde{\psi}(\xi)$ is said to be a
\emph{Stieltjes functions} if
%
\begin{equation}\label{eqstieltjes}
\tilde{\psi}(\xi) = \frac{\tilde{c}_1}{\xi} + \tilde{c_2} +
\frac{1}{\pi} \int_{0+}^\infty\frac{1}{\xi+ \zeta} \tilde{\mu
}(d\zeta) , \qquad \xi\in\C\setminus(-\infty, 0] ,
\end{equation}
for some $\tilde{c}_1, \tilde{c}_2 \ge0$ and some measure $\tilde
{\mu}$ on $(0, \infty)$ such that the integral $\int_0^\infty\min
(1, \zeta^{-1}) \tilde{\mu}(d\zeta)$ is finite. See \cite
{bibssv10} for a general account on complete Bernstein functions,
Stieltjes functions and related notions.

It is known that $\psi(\xi)$ is a CBF if and only if $\psi(\xi)$ is
nonnegative and increasing on $(0, \infty)$, holomorphic in $\C
\setminus(-\infty, 0]$, and $\im\psi(\xi) > 0$ when $\im\xi> 0$.
Furthermore, if $\psi(\xi)$ is a CBF, then $\xi/ \psi(\xi)$ is a
CBF, and $1 / \psi(\xi)$ and $\psi(\xi) / \xi$ are Stieltjes functions

The function $\tilde{\psi}(\xi)$ given by~\eqref{eqstieltjes} is
the Laplace transform of $\tilde{c}_2 \delta_0(dx) + (\tilde{c}_1 +
\laplace\tilde{\mu}(x))\,dx$ (\cite{bibssv10}, Theorem~2.2).
Furthermore, $\pi\tilde{c}_1 \delta_0(d\zeta) + \tilde{\mu
}(d\zeta)$ is the limit of measures $-\im(\tilde{\psi}(-\zeta+ i
\eps))\,d\zeta$ as $\eps\to0^+$ (\cite{bibssv10}, Corollary~6.3 and
Comments~6.12), so, in a sense, it is the boundary value of $\tilde
{\psi}$. Therefore, we use a shorthand notation $-\im(\tilde{\psi
}^+(-\zeta))\,d\zeta$ for $\tilde{\mu}(d\zeta)$. Furthermore, we
have $\tilde{c}_1 = \lim_{\xi\to0} (\xi\tilde{\psi}(\xi))$ and
$\tilde{c}_2 = \lim_{\xi\to\infty} \tilde{\psi}(\xi)$.

Following~\cite{bibmk10}, we define
%
\begin{equation}\label{eqdagger}
\psi^\dagger(\xi) = \exp\biggl( \frac{1}{\pi} \int_0^\infty \frac
{\xi\log\psi(\zeta^2)}{\xi^2 + \zeta^2} \, d\zeta \biggr) , \qquad
\re\xi> 0 ,
\end{equation}
for any function $\psi(\xi)$ such that $\min(1, \zeta^{-2}) \log
\psi(\zeta^2)$ is integrable in $\zeta> 0$. By a simple substitution,
%
\begin{equation}\label{eqdagger1}
\psi^\dagger(\xi) = \exp\biggl( \frac{1}{\pi} \int_0^\infty \frac
{\log\psi(\xi^2 \zeta^2)}{1 + \zeta^2} \, d\zeta \biggr) , \qquad \xi
> 0 .
\end{equation}
By~\cite{bibmk10}, Lemma~4, if $\psi(\xi)$ is a CBF, then also $\psi
^\dagger(\xi)$ is a CBF (this was independently proved in \cite
{bibksv10}, Proposition~2.4), and
%
\begin{equation}\label{eqduality}
\psi^\dagger(\xi) \psi^\dagger(-\xi) = \psi(-\xi^2) , \qquad
\xi
\in\C\setminus\R.
\end{equation}

\begin{proposition}
\label{propdaggerest}
If $\psi(\xi)$ is nonnegative on $(0, \infty)$, and both $\psi(\xi
)$ and $\xi/ \psi(\xi)$ are increasing on $(0, \infty)$, then
%
\begin{equation}\label{eqdaggerest}
e^{-2 \catalan/ \pi} \sqrt{\psi(\xi^2)} \le\psi^\dagger(\xi)
\le e^{2 \catalan/ \pi} \sqrt{\psi(\xi^2)} ,
\end{equation}
where $\catalan\approx0.916$ is the Catalan constant. Note that $e^{2
\catalan/ \pi} \le2$.

If, in addition, $\psi(\xi)$ is regularly varying at $\infty$, then
%
\begin{equation}\label{eqregvar}
\psi^\dagger(\xi) \sim\sqrt{\psi(\xi^2)}, \qquad \xi\to\infty.
\end{equation}
An analogous statement for $\xi\to0$ holds for $\psi(\xi)$
regularly varying at $0$.

In particular, \eqref{eqdaggerest} holds for any CBF. Likewise,~\eqref
{eqregvar} holds for any regularly varying CBF.
\end{proposition}

A result similar to~\eqref{eqdaggerest} was obtained independently
in~\cite{bibksv11}, Proposition~3.7, while~\eqref{eqregvar} for CBFs
was derived in~\cite{bibksv9}, Proposition 2.2.

\begin{pf}
By the assumptions, we have
%
\begin{equation}\label{eqpsiest}
\psi(\xi^2) \min(1, \zeta^2) \le\psi(\xi^2 \zeta^2) \le\psi
(\xi^2) \max(1, \zeta^2), \qquad \xi, \zeta> 0.
\end{equation}
It follows that
\begin{eqnarray*}
\psi^\dagger(\xi) & =& \exp\biggl( \frac{1}{\pi} \int_0^\infty \frac
{\log\psi(\xi^2 \zeta^2)}{1 + \zeta^2} \, d\zeta \biggr) \\ & \le&
\sqrt{\psi(\xi^2)} \exp\biggl( \frac{1}{\pi} \int_1^\infty\frac
{\log\zeta^2}{1 + \zeta^2}\,d\zeta \biggr) = e^{2 \catalan/ \pi} \sqrt
{\psi(\xi^2)} .
\end{eqnarray*}
The lower bound is obtained in a similar manner.

The second statement of the proposition is proved in a very similar
manner to Lemma~15 in~\cite{bibmk10}. Define an auxiliary function
$h(\xi, \zeta) = \psi(\xi^2 \zeta^2) / \psi(\xi^2)$. By~\eqref
{eqpsiest} we have $|\log h(\xi, \zeta)| \le2 |\log\zeta|$, $\xi,
\zeta> 0$. Since $\psi$ is regularly varying at~infinity, for some
$\alpha$, $\lim_{\xi\to\infty} h(\xi, \zeta) = \zeta^{2 \alpha
}$ for each $\zeta> 0$. Hence, by dominated convergence,
\[
\lim_{\xi\to\infty} \int_0^\infty\frac{\log h(\xi, \zeta)}{1 +
\zeta^2} \, d\zeta= \int_0^\infty\frac{\log\zeta^{2 \alpha
}}{1 + \zeta^2} \, d\zeta= 0 .
\]
It follows that
\[
\lim_{\xi\to\infty} \biggl( \int_0^\infty\frac{\log\psi(\xi^2 \zeta
^2)}{1 + \zeta^2} \, d\zeta- \frac{\pi}{2} \log\psi(\xi ^2) \biggr) =
0 ,
\]
and so finally $\lim_{\xi\to\infty} \psi^\dagger(\xi) / \sqrt
{\psi(\xi^2)} = 1$, as desired. Regular variation at $0$ is proved in
a similar way.
\end{pf}

As in~\cite{bibmk10}, for differentiable functions $\psi(\xi)$ with
positive derivative, we define
%
\begin{equation}\label{eqpsilambda}
\psi_\lambda(\xi) = \frac{1 - \xi/ \lambda^2}{1 - \psi(\xi) /
\psi(\lambda^2)} , \qquad\lambda> 0, \xi\in\C\setminus(-\infty
, 0) .
\end{equation}
This definition is extended continuously by $\psi_\lambda(\lambda^2)
= \psi(\lambda^2) / (\lambda^2 \psi'(\lambda^2))$. Note that if
$\psi(0) = 0$, then $\psi_\lambda(0) = 1$. For simplicity, we denote
$\psi_\lambda^\dagger(\xi) = (\psi_\lambda)^\dagger(\xi)$.
By~\cite{bibmk10}, Lemma~2, if $\psi(\xi)$ is a CBF, then $\psi
_\lambda(\xi)$ is a CBF for any $\lambda> 0$.

\subsection{Estimates for the Laplace transform}
\label{seclap}

This short section contains some rather standard estimates for the
inverse Laplace transform.

\begin{proposition}
\label{proplb}
Let $a > 0$, $c \ge1$. If $f$ is nonnegative and $f(x) \le c f(a) \max
(1, x/a)$ \mbox{($x > 0$)}, then for any $\xi> 0$,
\[
f(a) \ge\frac{\xi\laplace f(\xi)}{c (1 + (a \xi)^{-1} e^{-a \xi
})} .
\]
\end{proposition}

\begin{pf}
We have
\begin{eqnarray*}
\xi\laplace f(\xi) & =& \int_0^a \xi e^{-\xi x} f(x)\,dx + \int
_a^\infty\xi e^{-\xi x} f(x)\,dx \\
& \le& c f(a) \int_0^a \xi e^{-\xi x}\,dx + \frac{c f(a)}{a} \int
_a^\infty\xi x e^{-\xi x}\,dx \\
& =& c f(a) (1 - e^{-a \xi}) + \frac{c f(a)}{a \xi} (1 + a \xi)
e^{-a \xi} = c f(a) \bigl(1 + (a \xi)^{-1} e^{-a \xi}\bigr) ,
\end{eqnarray*}
as desired.
\end{pf}

\begin{proposition}
\label{propub}
If $f$ is nonnegative and increasing, then for $a, \xi> 0$,
\[
f(a) \le e^{a \xi} \xi\laplace f(\xi) .
\]
\end{proposition}

\begin{pf}
As before,
\begin{eqnarray*}
\xi\laplace f(\xi) & =& \int_0^a \xi e^{-\xi x} f(x)\,dx + \int
_a^\infty\xi e^{-\xi x} f(x)\,dx \\
& \ge& f(a) \int_a^\infty\xi e^{-\xi x}\,dx = f(a) e^{-a \xi} ,
\end{eqnarray*}
as claimed.\vadjust{\goodbreak}
\end{pf}

\begin{proposition}
\label{propub3}
If $f$ is nonnegative and decreasing, then for $a, \xi> 0$,
\[
f(a) \le\frac{\xi\laplace f(\xi)}{1 - e^{-a \xi}} .
\]
\end{proposition}

\begin{pf}
Again,
\begin{eqnarray*}
\xi\laplace f(\xi) & =& \int_0^a \xi e^{-\xi x} f(x)\,dx + \int
_a^\infty\xi e^{-\xi x} f(x)\,dx \\
& \ge& f(a) \int_0^a \xi e^{-\xi x}\,dx = f(a) (1 - e^{-a \xi}) ,
\end{eqnarray*}
as claimed.
\end{pf}

\section{Suprema of general L{\'e}vy processes}
\label{secsup}

We briefly recall the basic notions of the fluctuation theory for L{\'
e}vy processes. Let $L_t$ be the local time of the process $X_t$
reflected at its supremum $M_t$, and denote by $L^{-1}_s$ the
right-continuous inverse of $L_t$, the ascending ladder-time process
for $X_t$. This is a (possibly killed) subordinator, and $H_s =
X(L^{-1}_s) = M(L^{-1}_s)$ is another (possibly killed) subordinator,
called the ascending ladder-height process. The Laplace exponent of the
ascending ladder process, that is, the (possibly killed) bivariate
subordinator $(L^{-1}_s, H_s)$ ($s < L(\infty)$), is denoted by
$\kappa(z, \xi)$. By~\cite{bibb96}, Corollary~VI.10,
%
\begin{equation}\label{eqkappaFrulfull}
\kappa(z,\xi) = c\exp\biggl(\int_0^\infty\int_{[0,\infty
)}(e^{-t}-e^{-zt-\xi x})t^{-1}\pr(X_t\in dx )\,dt \biggr),
\end{equation}
where $c$ is a normalization constant of the local time. Since our
results are not affected by the choice of $c$, we assume that $c = 1$.
We note that $\kappa(z, 0)$ is a Bernstein function of $z$, and also
$z / \kappa(z, 0)$ is a Bernstein function (this follows from~\eqref
{eqkappaFrulfull} by Frullani's integral; see~\cite{bibb96},
formula (VI.3) for the case when $X_t$ is not a compound Poisson
process). For a more in-depth account of the fluctuation theory, we
refer the reader to~\cite{bibb96,bibd07,bibk06}. In general, there is
no closed-form formula for $\kappa(z, \xi)$. For a list of special
cases, see~\cite{bibkkp10} and the references therein. For a symmetric
process which is not a compound Poisson process, we have $\kappa(z,0)
= \sqrt{z}$.

As usual, $\tau_x$ denotes the first passage time through a barrier at
$x \ge0$ for~$X_t$ (or for $M_t$). Following~\cite{bibb96}, for $x, z
\ge0$, we define
\[
V^z(x) = \ex\biggl( \int_0^{\infty} \exp(-z L_s^{-1}) \ind_{[0,
x)}(H_s)\,ds \biggr) = \ex\biggl( \int_0^{\infty} e^{-z t} \ind_{[0, x)}(M_t)\,
dL_t \biggr) .
\]
For $z = 0$, we simply have $V^0(x) = \int_0^{\infty}\pr(H_s <
x)\,ds$, so that $V^0(x) = V(x)$ is the renewal function of the process
$H_s$, studied in more detail for symmetric L\'e{}vy processes in
Section~\ref{secv}. By~\cite{bibb96}, formula~(VI.8),
%
\begin{equation}\label{eqsuplaplaceta}
\int_0^\infty e^{-z t} \pr(M_t < x)\,dt = \frac{\kappa(z,0)
V^z(x)}{z} , \qquad x, z \ge0 .
\end{equation}
(Note that in~\cite{bibb96}, a weak inequality $M_t \le x$ is used in
the definition of $V^z(x)$.) Hence, for a symmetric process $X_t$ which
is not a compound Poisson process, we have
%
\begin{equation}\label{eqsuplaplacetb}
\int_0^\infty e^{-z t} \pr(M_t < x)\,dt = \frac{V^z(x)}{\sqrt{z}}
, \qquad x, z \ge0 .
\end{equation}
This is a partial inverse of the double Laplace transform in~\eqref
{eqbd}; however, there is no known explicit formula for $V^z(x)$. For a
different and, in a sense, more explicit partial inverse, see~\eqref
{eqsuplaplaceint} below.

By~\cite{bibb96}, Section~VI.4, the Laplace transform of $V^z(x)$ is
$1 / (\xi\kappa(z, \xi))$. Hence, when $X_t$ is symmetric and it is
not a compound Poisson process, the right-hand side of the
Baxter--Donsker formula~\eqref{eqbd} can be written as $\sqrt{z} / (z
\kappa(z, \xi))$; see~\cite{bibf74}, Corollary~9.7.

\begin{theorem}
\label{thsupest}
Let $X_t$ be a L{\'e}vy process, $M_t = \sup_{0 \le s \le t} X_s$ and
let $\kappa(z, \xi)$ be the bivariate Laplace exponent of its
ascending ladder process. Suppose that
%
\begin{equation}\label{eqKdef}
K(s) = \int_s^\infty\frac{\kappa(z, 0)}{z^2} \, dz < \infty,
\qquad
s > 0
\end{equation}
and that $\kappa(z, 0) / z$ is unbounded (near $0$). For $t, x > 0$,
we have
%
\begin{eqnarray}\label{eqsupesta}
\min\bigl( C_1, C_2(\kappa, t) \kappa(1 / t, 0) V(x) \bigr) & \le&\pr(M_t
< x) \nonumber
\\[-8pt]
\\[-8pt]
& \le&\min\biggl( 1, \frac{e}{e - 1} \kappa(1/t, 0) V(x) \biggr) .
\nonumber
\end{eqnarray}
Here
\[
C_1 = \frac{e - 1}{8 e^2} \quad \mbox{and} \quad C_2(\kappa, t) =
\frac
{z t}{2 e} ,
\]
where $z \in(0, 1/t)$ solves
\[
\frac{\kappa(z, 0)}{z} = \frac{4 e^2}{e - 1} K(1/t) .
\]
\end{theorem}

\begin{pf}
The upper bound in~\eqref{eqsupesta} is a direct consequence of~\eqref
{eqsuplaplaceta} and Proposition~\ref{propub3} with $\xi= 1/t$.

Following~\cite{bibb96}, Lemma~VI.21, we find a lower bound for
$V^z(x)$. We have
\begin{eqnarray*}
V(x) & =& \ex\biggl( \int_0^{\infty} \ind_{[0, x)}(M_t)\,dL_t \biggr) \\
& \le& e \ex\biggl( \int_0^{1/z} e^{-z t} \ind_{[0, x)}(M_t)\,dL_t \biggr) +
\ex\biggl( \int_{1/z}^\infty\ind_{[0, x)}(M_t)\,dL_t \biggr),
\end{eqnarray*}
which implies
%
\begin{equation}\label{eqVlower}
e V^z(x) \ge V(x) - \ex\biggl( \int_{1/z}^\infty\ind_{[0, x)}(M_t)\,dL_t \biggr).\vadjust{\goodbreak}
\end{equation}
Let $\sigma_z = \inf\{ t \ge1/z\dvtx X_t = M_t \} = L^{-1}(L_{1/z})$;
$\sigma_z$ is a stopping time. Since the support of the measure $dL_t$
is contained in the set $\{t \dvtx X_t = M_t\}$ of zeros of the reflected
process, we have
\begin{eqnarray*}
\ex\biggl( \int_{1/z}^\infty\ind_{[0, x)}(M_t)\,dL_t \biggr) & =& \ex\biggl( \int
_{\sigma_z}^\infty\ind_{[0, x)}(M_t)\,dL_t; M_{1/z} < x \biggr) \\
& \le&\ex\biggl( \int_{\sigma_z}^\infty\ind_{[0, x)}(M_t - M_{\sigma
_z})\,dL_t; M_{1/z} < x \biggr) .
\end{eqnarray*}
Next, observe that $M_{\sigma_z} = X_{\sigma_z}$, so that
\[
M_t - M_{\sigma_z} = \sup_{s \le t - \sigma_z} (X_{\sigma_z + s} -
X_{\sigma_z}) , \qquad t \ge\sigma_z .
\]
Hence,
\begin{eqnarray*}
&& \ex\biggl( \int_{1/z}^\infty\ind_{[0, x)}(M_t)\,dL_t \biggr) \\ && \qquad
\le\ex
\biggl( \int_{\sigma_z}^\infty\ind_{[0, x)}\Bigl( \sup_{s \le t-\sigma_z}
(X_{\sigma_z+s} - X_{\sigma_z}) \Bigr)\,dL_t; M_{1/z} < x \biggr) \\
&& \qquad = \ex\biggl( \int_0^\infty\ind_{[0, x)}\Bigl( \sup_{s \le u}
(X_{\sigma_z+s} - X_{\sigma_z}) \Bigr)\, d(L_{\sigma_z + u} - L_{\sigma
_z}); M_{1/z} < x \biggr) .
\end{eqnarray*}
Since $\sigma_z \ge1 / z$, by the strong Markov property,
\begin{eqnarray*}
\ex\biggl( \int_{1/z}^\infty\ind_{[0, x)}(M_t)\,dL_t \biggr) & \le&\pr
(M_{1/z} < x) \ex\biggl( \int_0^\infty\ind_{[0, x)}(M_u)\,dL_u \biggr) \\ &
=& \pr(M_{1/z} < x) V(x),
\end{eqnarray*}
which, by~\eqref{eqVlower}, yields
\[
V^z(x) \ge\frac{(1 - \pr(M_{1/z} < x)) V(x)}{e} = \frac{\pr
(M_{1/z} \ge x) V(x)}{e} .
\]
Let $k > 0$. By~\eqref{eqsuplaplaceta} and the already proved upper
bound of~\eqref{eqsupesta},
\begin{eqnarray*}
V^z(x) \kappa(z, 0) & =& z \int_0^{k/z} e^{-z t} \pr(M_t < x)\,dt + z
\int_{k/z}^\infty e^{-z t} \pr(M_t < x)\,dt \\
& \le&\frac{e}{e - 1} V(x) z \int_0^{k/z} e^{-z t} \kappa(1/t, 0)\,dt
+ \pr(M_{k/z} < x) .
\end{eqnarray*}
The last two estimates give
%
\begin{eqnarray}\label{eqprest1}\qquad
\pr(M_{k/z} < x) & \ge&\frac{\kappa(z, 0) \pr(M_{1/z} \ge x)
V(x)}{e} - \frac{e}{e - 1} V(x) z \int_0^{k/z} \kappa(1/t, 0)\,
dt\nonumber
\\[-8pt]
\\[-8pt]
& =& \frac{V(x) \kappa(z,0)}{e} \biggl( \pr(M_{1/z} \ge x) - \frac
{e^2}{e - 1} \frac{z K(z/k)}{\kappa(z, 0)} \biggr).
\nonumber
\end{eqnarray}
Fix $\eps\in(0, 1)$ (later we choose $\eps= 1/4$). Note that the
function $\kappa(z, 0) / z$ is continuous, decreasing and unbounded.
Hence, it maps the interval $(0, 1/t]$ onto the interval $[t\kappa
(1/t, 0), \infty)$. Furthermore, $\kappa(z, 0)$ is increasing, so
that $K(z) \ge\kappa(z, 0) / z$. In particular, $\frac{e^2}{\eps(e
- 1)} K(1/t)>K(1/t) \ge t\kappa(1/t, 0) $. It follows that we can
choose $z = z(t)<1 / t$ such that
\[
\frac{\kappa(z, 0)}{z} = \frac{e^2}{\eps(e - 1)} K(1/t).
\]
Setting $k = z t < 1$, the above equality can be rewritten as
%
\begin{equation}\label{eqprest2}
\frac{e^2}{e - 1} \frac{z K(z / k)}{\kappa(z, 0)} = \eps.
\end{equation}
Suppose now that $V(x) \kappa(z, 0) \le\eps(e - 1) / e$. Then, by
the upper bound of~\eqref{eqsupesta}, we have $\pr(M_{1/z} \ge x) = 1
- \pr(M_{1/z} < x) \ge1 - \eps$. This, \eqref{eqprest1} and~\eqref
{eqprest2} give
\[
\pr(M_t < x) = \pr(M_{k/z} < x) \ge\frac{V(x) \kappa(z,0)}{e} (1
- 2 \eps) .
\]
This estimate holds for $t \ge t_0$, where $V(x) \kappa(z(t_0), 0) =
\eps(e - 1) / e$ [here we use continuity of $\kappa(z(t), 0)$ as a
function of $t$]. Hence, by monotonicity of $\pr(M_t < x)$ in $t$,
\[
\pr(M_t < x) \ge\min\biggl( \frac{\eps(1 - 2 \eps) (e - 1)}{e^2}, \frac
{(1 - 2 \eps) V(x) \kappa(z,0)}{e} \biggr) .
\]
The lower bound in~\eqref{eqsupesta} follows by taking $\eps= 1/4$
and using the inequality $\kappa(z, 0) = \kappa(k/t, 0) \ge k \kappa
(1/t, 0)$.
\end{pf}

To formulate the next result we define the following \emph{upper
scaling conditions}:
%
\begin{eqnarray}
\label{equscA}
\mbox{for some $\varrho\in(0, 1)$ and $c > 0$,} \qquad  \frac
{\kappa(z_2, 0)}{\kappa(z_1, 0)} &\le& c \frac{z_2^\varrho
}{z_1^\varrho}  \qquad  \qquad  \nonumber
\\[-8pt]
\\[-8pt]
\eqntext{\mbox{when $0 < z_1 < z_2 < 1$,}} \\
\label{equscB}
\mbox{for some $\varrho\in(0, 1)$ and $c > 0$,} \qquad  \frac
{\kappa(z_2, 0)}{\kappa(z_1, 0)} &\le& c \frac{z_2^\varrho
}{z_1^\varrho} \qquad  \qquad  \nonumber
\\[-8pt]
\\[-8pt]
\eqntext{\mbox{when $1 < z_1 < z_2$.}}
\end{eqnarray}
Observe that condition \eqref{equscB} implies that for any $z^* > 0$,
there is $c^*$ such that
%
\begin{equation}\label{equscB1}
\frac{\kappa(z_2, 0)}{\kappa(z_1, 0)} \le c^* \frac{z_2^\varrho
}{z_1^\varrho}  \qquad \mbox{when $z^* < z_1 < z_2$.}
\end{equation}

\begin{corollary}
\label{corkapparegular1}
Let $X_t$ be a L{\'e}vy process, $M_t = \sup_{0 \le s \le t} X_s$ and
let $\kappa(z, \xi)$ be the bivariate Laplace\vadjust{\goodbreak} exponent of its
ascending ladder process. If $\kappa(z, 0)$ satisfies condition~\eqref
{equscA} with $0 < \ro< 1$ and the integral $\int_1^\infty\kappa(z,
0) z^{-2}\,dz$ is finite, then
%
\begin{equation}\label{eqsupest}
  C(\kappa) \min\bigl( 1 ,\kappa(1 / t, 0) V(x) \bigr) \le\pr(M_t < x) \le
\min\bigl( 1, 2 \kappa(1 / t, 0) V(x) \bigr) ,\hspace*{-35pt}
\end{equation}
for every $x > 0$ and $t \ge1$. If $\kappa(z, 0)$ satisfies~\eqref
{equscB} with $0 < \ro< 1$ and $\lim_{z \to0} z / \kappa(z ,0) =
0$, then~\eqref{eqsupest} holds for $x > 0$ and $t \le1$.

In particular, if $\kappa(z, 0)$ satisfies both~\eqref{equscA}
and~\eqref{equscB}, that is, there are $c > 0$ and $\ro\in(0, 1)$
such that $\kappa(\lambda z, 0) \le c \lambda^\ro\kappa(z, 0)$ for
$\lambda\ge1$ and $z > 0$, then~\eqref{eqsupest} is true for every
$x > 0$ and $t > 0$.
\end{corollary}

\begin{pf}
We begin with the first part of the statement. By condition~\eqref{equscA},
\[
\kappa(z, 0) \le c_1(\kappa) \biggl( \frac z s \biggr)^\ro\kappa(s, 0),
\qquad s \le z \le1 .
\]
In particular, $\kappa(s, 0) / s$ is unbounded. Furthermore, using
also finiteness of the integral $\int_1^\infty\kappa(z, 0) z^{-2}\,
dz$, we obtain
%
\begin{equation}\label{eqKest}
K(s) \le c_2(\kappa) \frac{\kappa(s, 0)}{s} , \qquad s \le1.
\end{equation}
This implies that the assumptions of Theorem~\ref{thsupest} are satisfied.

Let $t \ge1$ and define $z = z(t) \in(0, 1/t)$ as in Theorem~\ref
{thsupest}. By condition~\eqref{equscA} we have
\[
\frac{\kappa(1/t, 0)}{\kappa(z, 0)} \le\frac{c_3(\kappa)}{(z
t)^{\ro}} .
\]
By definition of $z$ and~\eqref{eqKest} (with $s = 1 / t$), we have
\[
\frac{1}{z} = \frac{4 e^2}{e - 1} \frac{K(1/t)}{\kappa(z, 0)} \le
\frac{4 e^2 c_2(\kappa) c_3(\kappa)}{e - 1} \frac{t}{(t z)^{ \ro}} ,
\]
which gives $z t \ge c_4(\kappa)$. Hence, the constant $C_2$ in
Theorem~\ref{thsupest} satisfies $C_2 = z t / (2 e) \ge c_4(\kappa) /
(2 e)$. This ends the proof of the first part.

The second part can be justified in a similar way, since
condition~\eqref{equscB} implies that
\[
K(s) \le c_5(\kappa) \frac{\kappa(s, 0)}{s} , \qquad s \ge1 .
\]
Moreover, for $t < 1$ and $z = z(t)$ selected according to Theorem~\ref
{thsupest} we have $z(1) \le z(t) < 1/t$. Applying \eqref{equscB}
[with $z^* = z(1)$], we obtain
\[
\frac{\kappa(1/t, 0)}{\kappa(z, 0)} \le\frac{c_6(\kappa)}{(z
t)^{\ro}} , \qquad z \le\frac{1}{t} .
\]
Finally, the last statement is a direct consequence of the previous ones.
\end{pf}

\begin{remark}
$\!\!$Due to Potter's theorem (\cite{bibbgt87}, Theorem~1.5.6)
condition~\eqref{equscA} is implied by regular variation of $\kappa
(z, 0)$ at zero with index $0 < \ro^* <1$.\vadjust{\goodbreak} Likewise, condition~\eqref
{equscB} is implied by regular variation of $\kappa(z, 0)$ at $\infty
$ with index $0 < \ro^* <1$.

In the second part of the above corollary the assumption $\lim_{z \to
0} z / \kappa(z ,0) = 0$ can be removed at the expence that the lower
bound holds for $t \le t_0$, where $t_0 =t_0(\kappa)$ is sufficiently
small. This is due to the fact that since $\lim_{t \searrow0} K(1/t)
= 0$, $z = z(t)$ in Theorem~\ref{thsupest} is well defined for $t$
small enough.
\end{remark}

By the results of~\cite{bibb96}, Theorem VI.14 and~\cite{bibbd97},
the regular variation of order $\ro\in(0, 1)$ of $\kappa(z, 0)$ at
$0$ or at $\infty$ is equivalent to the existence of the limit of $\pr
(X_t \ge0)$ as $t \to\infty$ or $t \to0^+$, respectively. Hence,
Corollary~\ref{corkapparegular1} implies the following result.

\begin{corollary}
\label{corkapparegular2}
Let $X_t$ be a L{\'e}vy process and $M_t = \sup_{0 \le s \le t} X_s$. If
\[
\lim_{t \to\infty} \pr(X_t \ge0) \in(0,1) \quad \mbox{and}
\quad
\limsup_{t \to0^+} \pr(X_t \ge0) < 1 ,
\]
then~\eqref{eqsupest} holds for $x > 0$ and $t \ge1$. If
\[
\lim_{t \to0^+} \pr(X_t \ge0) \in(0,1) \quad \mbox{and} \quad
\limsup
_{t \to\infty} \pr(X_t \ge0) < 1 ,
\]
then~\eqref{eqsupest} is true for $x>0$ and $t \le1$. Finally, if
\[
\lim_{t \to\infty} \pr(X_t \ge0) \in(0,1) \quad \mbox{and}
\quad \lim
_{t \to0^+} \pr(X_t \ge0) \in(0,1) ,
\]
then~\eqref{eqsupest} holds for every $x>0$ and $t>0$.
\end{corollary}

\begin{pf}
We only need to verify that $\kappa(z, 0) / z^2$ is integrable at
infinity, and that $\lim_{z \to0^+} (z / \kappa(z, 0)) = 0$. In each
of the cases, there is $\eps> 0$ such that $\pr(X_t \ge0) \le1 -
\eps$ for all $t > 0$. Therefore, by~\eqref{eqkappaFrulfull} and the
Frullani integral, $\kappa(z, 0) \le z^{1 - \eps}$ for $z \ge1$, and
$\kappa(z, 0) \ge z^{1 - \eps}$ when $0 < z < 1$. The result follows.
\end{pf}

\begin{remark}
The uniform estimates of Corollary~\ref{corkapparegular2} complement
the existing results from~\cite{bibgn86} about the asymptotic behavior
of $\pr(M_t < x)$, where it was shown that
\[
\lim_{t \to\infty} \frac{\sqrt{\pi}}{\kappa(1/t,0)} \pr(M_t <
x) = V(x) ,
\]
under the assumption that $\kappa(z,0)$ is regularly varying at zero
with index $\ro\in(0, 1)$.
\end{remark}

\section{Suprema of symmetric L{\'e}vy processes}
\label{secv}

In this section we assume that $X_t$ is a \emph{symmetric} L{\'e}vy
process with L{\'e}vy--Khintchin exponent $\Psi(\xi)$.
In a rather general setting, we can invert the Laplace transform in
time variable in~\eqref{eqbd}.

\begin{theorem}
\label{thsuplaplace}
Suppose that $X_t$ is a symmetric L{\'e}vy process with L{\'
e}vy--Khintchin exponent\vadjust{\goodbreak} $\Psi(\xi)$. Suppose that $\Psi(\xi)$ is
increasing in $\xi> 0$. If $M_t = \sup_{0 \le s \le t} X_s$, then
%
\begin{eqnarray}\label{eqsuplaplace}
\ex e^{-\xi M_t} & =& \frac{1}{\pi} \int_0^\infty\frac{\xi\Psi
'(\lambda)}{(\lambda^2 + \xi^2) \sqrt{\Psi(\lambda)}} \nonumber
\\[-8pt]
\\[-8pt]
&&\hphantom{\frac{1}{\pi} \int_0^\infty}{} \times\exp\biggl( \frac{1}{\pi} \int_0^\infty\frac{\xi\log
\fracf{\lambda^2 - \zeta^2}{\Psi(\lambda) - \Psi(\zeta)}}{\xi^2
+ \zeta^2} \, d\zeta \biggr) e^{-t \Psi(\lambda)}\,d\lambda.
\nonumber
\end{eqnarray}
\end{theorem}

Since $\pr(M_t < x) = \pr(\tau_x > t)$, the following integrated
form of~\eqref{eqsuplaplace} is sometimes more convenient.

\begin{corollary}
\label{corsuplaplace}
With the notation and assumptions of Theorem~\ref{thsuplaplace},
%
\begin{eqnarray}\label{eqsuplaplaceint}\quad
&&\int_0^\infty e^{-\xi x} \pr(\tau_x > t)\,dx\nonumber\\
 && \qquad  = \frac{\ex e^{-\xi
M_t}}{\xi}\nonumber
\\[-8pt]
\\[-8pt]
&& \qquad = \frac{1}{\pi} \int_0^\infty\frac{\Psi'(\lambda
)}{(\lambda^2 + \xi^2) \sqrt{\Psi(\lambda)}}
\nonumber \\
&& \qquad  \quad \hphantom{\frac{1}{\pi} \int_0^\infty}{} \times\exp\biggl( \frac{1}{\pi} \int_0^\infty\frac{\xi\log
(\fracf{\lambda^2 - \zeta^2}{\Psi(\lambda) - \Psi(\zeta)})}{\xi^2
+ \zeta^2} \, d\zeta \biggr) e^{-t \Psi(\lambda)}\,d\lambda.
\nonumber
\end{eqnarray}
\end{corollary}

\begin{pf*}{Proof of Theorem~\ref{thsuplaplace}}
Let $\psi(\xi) = \Psi(\sqrt{\xi})$ for $\xi> 0$. For any $z \in
\C\setminus(-\infty, 0]$ and $\xi> 0$, we define [see~\eqref{eqbd}
and~\eqref{eqdagger}]
\begin{eqnarray*}
\ph(\xi, z) & =& \sqrt{z} \exp\biggl( -\frac{1}{\pi} \int_0^\infty
\frac{\xi\log(z + \Psi(\zeta))}{\xi^2 + \zeta^2} \, d\zeta \biggr) \\
& =& \exp\biggl( -\frac{1}{\pi} \int_0^\infty\frac{\xi\log(1 + \fraca
{\psi(\zeta^2)}{z})}{\xi^2 + \zeta^2}\,d\zeta \biggr) .
\end{eqnarray*}
For any $\xi> 0$, the function $\ph(\xi, z)$ is positive and
increasing in $z \in(0, \infty)$. As $z \to0$ or $z \to\infty$,
$\ph(\xi, z)$ converges to $0$ and $1$, respectively. Furthermore, if
$\im z > 0$, then $\arg(1 + \psi(\zeta^2) / z) \in(-\pi, 0)$ for
all $\zeta> 0$, and therefore
\[
\arg\ph(\xi, z) = -\frac{1}{\pi} \int_0^\infty\frac{\xi\arg
(1 + \fraca{\psi(\zeta^2)}{z})}{\xi^2 + \zeta^2} \, d\zeta\in(0,
\pi/2) .
\]
Hence, for any $\xi> 0$, $\ph(\xi, z)$ [and even $(\ph(\xi,
z))^2$] is a complete Bernstein function of $z$. Note that the
continuous boundary limit $\ph^+(\xi, -z)$ exists for $z > 0$: if $z
= \psi(\lambda^2)$, or $\lambda= \sqrt{\psi^{-1}(z)}$, then
\begin{eqnarray*}
\ph^+(\xi, -z) & =& \exp\biggl( -\frac{1}{\pi} \int_0^\infty\frac
{\xi\log^- (1 - \fraca{\psi(\zeta^2)}{\psi(\lambda^2)})}{\xi^2 +
\zeta^2} \, d\zeta \biggr) \\
& =& \exp\biggl( -\frac{1}{\pi} \int_0^\infty\frac{\xi\log|1 - \fraca
{\psi(\zeta^2)}{\psi(\lambda^2)}|}{\xi^2 + \zeta^2} \, d\zeta+ i
\int_\lambda^\infty\frac{\xi}{\xi^2 + \zeta^2}\,d\zeta \biggr) \\
& =& \exp\biggl( \frac{1}{\pi} \int_0^\infty\frac{\xi(\log\psi
_\lambda(\zeta^2) - \log|1 - \fraca{\zeta^2}{\lambda^2}|)}{\xi^2
+ \zeta^2} \, d\zeta+ i \arctan\frac{\xi}{\lambda} \biggr) ;
\end{eqnarray*}
see~\eqref{eqpsilambda} for the notation. Here $\log^-$ denotes the
boundary limit on $(-\infty, 0)$ approached from below, $\log
^-(-\zeta) = -i \pi/2 + \log\zeta$ for $\zeta> 0$. The function
$\log|1 - \zeta^2 / \lambda^2|$ is harmonic in the upper half-plane
$\im\zeta> 0$, so that
\[
\frac{1}{\pi} \int_0^\infty\frac{\xi\log|1 - \fraca{\zeta
^2}{\lambda^2}|}{\xi^2 + \zeta^2} \, d\zeta= \frac{1}{2} \log
\biggl( 1 + \frac{\xi^2}{\lambda^2} \biggr) .
\]
Furthermore, $\exp(i \arctan(\xi/\lambda)) = (\lambda+ i \xi) /
\sqrt{\lambda^2 + \xi^2}$. Therefore, with $z = \psi(\lambda^2)$,
%
\begin{eqnarray}\label{eqphiplus}
\ph^+(\xi, -z) & =& \frac{\lambda(\lambda+ i \xi)}{\lambda^2 +
\xi^2} \exp\biggl( \frac{1}{\pi} \int_0^\infty\frac{\xi\log\psi
_\lambda(\zeta^2)}{\xi^2 + \zeta^2}\,d\zeta \biggr) \nonumber
\\[-8pt]
\\[-8pt] & =& \frac{\lambda
(\lambda+ i \xi) \psi_\lambda^\dagger(\xi)}{\lambda^2 + \xi^2} ;
\nonumber
\end{eqnarray}
see~\eqref{eqdagger} for the notation. Note that if $\psi(\xi)$ is
bounded on $(0, \infty)$ and $z \ge\sup_{\xi> 0} \psi(\xi)$, then
$\ph^+(\xi, -z)$ is real.

By~\eqref{eqbd}, $\ph(\xi, z) / z$ is the double Laplace transform
of the distribution of $M_t$. But for all $\xi> 0$, $\ph(\xi, z) /
z$ is a Stieltjes function of $z$. Therefore, by~\eqref{eqstieltjes},
\begin{eqnarray*}
\frac{\ph(\xi, z)}{z} & =& \frac{1}{\pi} \int_0^\infty\im\frac
{\ph^+(\xi, -\zeta)}{\zeta} \frac{1}{z + \zeta} \, d\zeta\\
& =& \frac{1}{\pi} \int_0^\infty2 \lambda\psi'(\lambda^2) \im
\frac{\ph^+(\xi, -\psi(\lambda^2))}{\psi(\lambda^2)} \frac{1}{z
+ \psi(\lambda^2)} \, d\lambda\\
& =& \frac{2}{\pi} \int_0^\infty\frac{\lambda\psi'(\lambda
^2)}{\psi(\lambda^2)} \frac{\lambda\xi\psi_\lambda^\dagger(\xi
)}{\lambda^2 + \xi^2} \frac{1}{z + \psi(\lambda^2)}\,d\lambda.
\end{eqnarray*}
Note that the second equality holds true also when $\psi(\xi)$ is
bounded. Since $1 / (z + \psi(\lambda^2)) = \int_0^\infty e^{-t \psi
(\lambda^2)} e^{-z t}\,dt$, we have
\[
\frac{\ph(\xi, z)}{z} = \int_0^\infty\biggl( \frac{2}{\pi} \int
_0^\infty\frac{\lambda\psi'(\lambda^2)}{\psi(\lambda^2)} \frac
{\lambda\xi\psi_\lambda^\dagger(\xi)}{\lambda^2 + \xi^2} e^{-t
\psi(\lambda^2)}\,d\lambda \biggr) e^{-z t}\,dt .
\]
The theorem follows by the uniqueness of the Laplace transform.
\end{pf*}

Let $V(x) = V^0(x)$ be the renewal function for the ascending
ladder-height process $H_s$ corresponding to $X_t$; see Section~\ref
{secsup} for the definition. When $X_t$ satisfies the \emph{absolute
continuity condition}\vadjust{\goodbreak} [e.g., if $1 / (1 + \Psi(\xi))$ is integrable
in $\xi$], then $V(x)$ is the (unique up to a multiplicative constant)
increasing harmonic function for $X_t$ on $(0, \infty)$, and $V'(x)$
is the decreasing harmonic function for $X_t$ on $(0, \infty)$;
cf.~\cite{bibs80}. It is known (\cite{bibb96}, formula~(VI.6)) that
for $\xi> 0$,
\[
\laplace V(\xi) = \frac{1}{\xi\kappa(0, \xi)} ,
\]
Moreover, if $X_t$ is not a compound Poisson process, then by \cite
{bibf74}, Corollary~9.7,
\[
\kappa(0, \xi) = \exp\biggl( \frac{1}{\pi} \int_0^\infty\frac {\xi
\log\Psi(\zeta)}{\xi^2 + \zeta^2} \, d\zeta \biggr) = \psi^\dagger
(\xi) ,
\]
where $\Psi(\xi) = \psi(\xi^2)$; see~\eqref{eqdagger} for the
notation. Clearly, we have $\laplace V'(\xi) = \xi\laplace V(\xi) =
1 / \psi^\dagger(\xi)$; here $V'$ is the distributional derivative
of $V$ on $[0, \infty)$. We remark that when $X_t$ is a compound
Poisson process, then, also by~\cite{bibf74}, Corollary~9.7,
%
\begin{equation}\label{eqcpp}
 \quad \kappa(0, \xi) = c \psi^\dagger(\xi) \qquad \mbox{with } c =
\exp\biggl( -\frac{1}{2} \int_0^\infty\frac{1 - e^{-t}}{t} \pr(X_t =
0)\,dt \biggr) .
\end{equation}
For simplicity, we state the next three results only for the case when
$X_t$ is not a compound Poisson process. However, extensions for
compound Poisson processes are straightforward due to~\eqref{eqcpp}.

As an immediate consequence of Proposition~\ref{propdaggerest} and
Karamata's Tauberian theorem (\cite{bibbgt87}, Theorem~1.7.1), we
obtain the following result, which in the case of complete Bernstein
functions was derived in Proposition 2.7 of~\cite{bibksv9}.

\begin{proposition}
\label{propvregular}
Let $\Psi(\xi)$ be the L{\'e}vy--Khintchin exponent of a symmetric
L{\'e}vy process $X_t$, which is not a compound Poisson process, and
suppose that both $\Psi(\xi)$ and $\xi^2 / \Psi(\xi)$ are
increasing in $\xi> 0$. If $\Psi(\xi)$ is regularly varying at
$\infty$, then $V$ is regularly varying at $0$ and $\Gamma(1 + \alpha
) V(x) \sim1 / \sqrt{\Psi(1/x)}$ as $x \to0$. Similarly, if $\Psi
(\xi)$ is regularly varying at $0$, then $\Gamma(1 + \alpha) V(x)
\sim1 /\sqrt{\Psi(1/x)}$ as $x \to\infty$.
\end{proposition}

Another consequence of Proposition~\ref{propdaggerest} is a uniform
estimate of the renewal function; see also Proposition~3.9 of~\cite{bibksv11}.

\begin{theorem}
\label{thvest}
Let $\Psi(\xi)$ be the L{\'e}vy--Khintchin exponent of a symmetric
L{\'e}vy process $X_t$, which is not a compound Poisson process, and
suppose that both $\Psi(\xi)$ and $\xi^2 / \Psi(\xi)$ are
increasing in $\xi> 0$. Then
\[
\frac{1}{5} \frac{1}{\sqrt{\Psi(1/x)}} \le V(x) \le5 \frac
{1}{\sqrt{\Psi(1/x)}} , \qquad x > 0 .
\]
\end{theorem}

\begin{pf}
Let $\psi(\xi) = \Psi(\sqrt{\xi})$ for $\xi> 0$. By
Proposition~\ref{propdaggerest}, we obtain $e^{-2 \catalan/ \pi} /
\sqrt{\xi^2 \psi(\xi^2)} \le\laplace{V}(\xi) \le e^{2 \catalan/
\pi} / \sqrt{\xi^2 \psi(\xi^2)}$, $\xi> 0$.\vadjust{\goodbreak} Since $V$ is
increasing, Proposition~\ref{propub} gives
\[
V(x) \le\frac{e \laplace V(1/x)}{x} \le\frac{e^{1 + 2 \catalan/
\pi}}{\sqrt{\psi(1/x^2)}} \le\frac{5}{\sqrt{\psi(1/x^2)}} .
\]
Furthermore, using subadditivity and monotonicity of $V$ (see \cite
{bibb96}, Section~III.1), for $x = k a + r$ ($k \ge0$, $r \in[0, a)$)
we obtain $V(x) \le k V(a) + V(r) \le(k + 1) V(a)$. It follows that
$V(x) \le2 V(a) \max(1, x / a)$ for all $a, x > 0$, and so, by
Proposition~\ref{proplb},
\[
V(x) \ge\frac{\laplace V(1/x)}{2 x (1 + e^{-1})} \ge\frac{1}{2 (1
+ e^{-1}) e^{2 \catalan/ \pi} \sqrt{\psi(1/x^2)}} \ge\frac{1}{5
\sqrt{\psi(1/x^2)}} ,
\]
as desired.
\end{pf}

We remark that when $V$ is a concave function on $(0, \infty)$ (e.g.,
when $\psi$ is a complete Bernstein function, see below), then clearly
$V(x) \le\max(1, x/a) V(a)$, so that the lower bound in Theorem~\ref
{thvest} holds with constant $2/5$ instead of~$1/5$.

If $\psi(\xi)$ is a complete Bernstein function [CBF, see~\eqref
{eqcbf}], then $\psi^\dagger(\xi)$ and $\xi/ \psi^\dagger(\xi)$
are CBFs, and hence $1 / \psi^\dagger(\xi)$ is a Stieltjes function;
see~\eqref{eqstieltjes}. Therefore, $V'(x)$ is a completely monotone
function on $(0, \infty)$, and $V(x)$ is a Bernstein function;
see~\cite{bibssv10} for the relation between completely monotone,
Bernstein, complete Bernstein and Stieltjes functions. More precisely,
we have the following result.

\begin{proposition}
\label{propv}
Let $\Psi(\xi)$ be the L{\'e}vy--Khintchin exponent of a symmetric
L{\'e}vy process $X_t$, which is not a compound Poisson process, and
suppose that $\Psi(\xi) = \psi(\xi^2)$ for a complete Bernstein
function $\psi$. Then $V$ is a Bernstein function, and
%
\begin{eqnarray}  \quad \qquad
\label{eqv}
V(x) & =& b x + \frac{1}{\pi} \int_{0^+}^\infty\im\biggl( -\frac
{1}{\psi^+(-\xi^2)} \biggr) \frac{\psi^\dagger(\xi)}{\xi} (1 - e^{-x
\xi})\,d\xi, \qquad x > 0 , \\
\label{eqvprime}
  \quad \qquad
V'(x) & =& b + \frac{1}{\pi} \int_{0^+}^\infty\im\biggl( -\frac
{1}{\psi^+(-\xi^2)} \biggr) \psi^\dagger(\xi) e^{-x \xi}\,d\xi,
\qquad x > 0 ,
\end{eqnarray}
where $b = \lim_{\xi\to0^+} (\xi/ \sqrt{\psi(\xi^2)})$.
\end{proposition}

As explained after formula~\eqref{eqstieltjes}, the expression $\im
(-1 / \psi^+(-\xi^2))\,d\xi$ in~\eqref{eqv} and~\eqref{eqvprime}
should be understood in the distributional sense, as a weak limit of
measures $\im(-1 / \psi(-\xi^2 + i \eps))\,d\xi$ on $(0, \infty)$
as $\eps\to0^+$. The measure $\im(-1 /\break \psi^+(-\xi))\,d\xi$ has an
atom of mass $\pi b$ at $0$, and this atom is not included in the
integrals from $0^+$ to $\infty$ in~\eqref{eqv} and~\eqref{eqvprime}.

\begin{pf}
Since $1 / \psi^\dagger(\xi)$ is a Stieltjes function, it has the
form~\eqref{eqstieltjes},
\[
\laplace V'(\xi) = \frac{1}{\psi^\dagger(\xi)} = a + \frac
{b}{\xi} + \frac{1}{\pi} \int_{0+}^\infty\frac{1}{\xi+ \zeta}
\tilde{\mu}(d\zeta) , \qquad \xi\in\C\setminus(-\infty, 0] ,\vadjust{\goodbreak}
\]
where, using~\eqref{eqduality},
\[
\tilde{\mu}(d\xi) = -\im\biggl( \frac{1}{(\psi^\dagger)^+(-\xi )} \biggr)\,
d\xi= -\im\biggl( \frac{\psi^\dagger(\xi)}{\psi^+(-\xi^2)} \biggr)\,d\xi
\]
and
\[
a = \lim_{\xi\to\infty} \frac{1}{\psi^\dagger(\xi)}
, \qquad b = \lim_{\xi\to0^+} \frac{\xi}{\psi^\dagger(\xi)}
.
\]
Using Proposition~\ref{propdaggerest}, we can express $a$ and $b$ in
terms of $\psi$. Since $\psi$ is unbounded, also $\psi^\dagger$ is
unbounded [by~\eqref{eqdaggerest}], and so in fact $a = 0$. In a
similar way, if $\xi/ \psi(\xi)$ converges to $0$ as $\xi\to0^+$,
then~\eqref{eqdaggerest} gives $\xi/ \psi^\dagger(\xi) \to0$, so
that $b = 0$. When the limit of $\xi/ \psi(\xi)$ is positive [since
$\xi/ \psi(\xi)$ is a CBF, the limit always exists], then $\psi$ is
regularly varying at $0$, and so $b = \lim_{\xi\to0^+} (\xi/ \sqrt
{\psi(\xi^2)})$, as desired. By the uniqueness of the Laplace transform,
\[
V'(x) = b + \frac{1}{\pi} \int_{0^+}^\infty e^{-x \xi} \tilde
{\mu}(d\xi) , \qquad x > 0 .
\]
The result follows by integration in $x$.
\end{pf}

Note that for a compound Poisson process, we have $a > 0$, so there is
an extra positive constant in~\eqref{eqv}.

As a combination of Theorem~\ref{thsupest} and Theorem~\ref{thvest},
we obtain the following result.

\begin{theorem}
\label{thestsym}
Let $\Psi(\xi)$ be the L{\'e}vy--Khintchin exponent of a symmetric
L{\'e}vy process $X_t$. Suppose that both $\Psi(\xi)$ and $\xi^2 /
\Psi(\xi)$ are increasing in $\xi> 0$. If $M_t = \sup_{0 \le s \le
t} X_s$, then for all $t, x > 0$,
\[
\frac{1}{100} \min\biggl( 1, \frac{1}{200 \sqrt{t \Psi(1/x)}} \biggr)
\le\pr(M_t < x) \le\min\biggl( 1, \frac{10}{\sqrt{t \Psi(1/x)}} \biggr) .
\]
\end{theorem}

\begin{pf}
When $X_t$ is not a compound Poisson process, then the result follows
from Theorems~\ref{thsupest} and~\ref{thvest}, and from $\kappa(z,
0) = \sqrt{z}$. Suppose that $X_t$ is a compound Poisson process. For
$\eps> 0$ consider $X^\eps_t = \eps B_t + X_t$, where the Brownian
motion $B_t$ is independent of $X_t$. Then the L{\'e}vy--Khintchin
exponent of $X^\eps_t$ equals to $\Psi_\eps(\xi)=(\eps\xi)^2 +
\Psi(\xi)$. It is easy to check that $\xi^2 / \Psi_\eps(\xi)$ is
increasing. Moreover, $M^\eps_t$ converges in distribution to $M_t$ as
$\eps\to0$. The result follows by the continuity of $\Psi(\xi)$.
\end{pf}

\begin{remark}
Clearly, the condition ``$\Psi(\xi)$ and $\xi^2 / \Psi(\xi)$ are
increasing in $\xi> 0$,'' in Theorem~\ref{thvest}, Proposition~\ref
{propvregular} and Theorem~\ref{thestsym}, can be replaced with
%
\begin{equation}\label{eqaltcond}
0 < \Psi'(\xi) < \frac{2 \Psi(\xi)}{\xi} , \qquad \xi> 0 .
\end{equation}
If $\Psi(\xi) = \psi(\xi^2)$, then~\eqref{eqaltcond} reads
%
\begin{equation}\label{eqaltcond2}
0 < \psi'(\xi) < \frac{\psi(\xi)}{\xi} , \qquad \xi> 0 .
\end{equation}
Using the standard representation of Bernstein functions, it is easy to
check that any Bernstein function
$\psi(\xi)$ (not necessarily a complete one) satisfies~\eqref
{eqaltcond2}. Hence, Theorem~\ref{thestsym}
applies to any \emph{subordinate Brownian motion}: a process $X_t =
B_{\eta_t}$, where $B(s)$ is the standard
Brownian motion [with $\ex(B_s) = 0$ and $\var(B_s) = 2 s$], $\eta
_t$ is a subordinator
[with $\ex(e^{-\xi\eta_t}) = e^{-t \psi(\xi)}$], and $B_s$ and
$\eta_t$ are independent processes.
\end{remark}

\section*{Acknowledgments}
We are deeply indebted to Lo{\"i}c Chaumont for numerous discussions on
the subject of the article and many valuable suggestions.
We thank Tomasz Grzywny for pointing out errors in the preliminary
version of the article. We also thank the anonymous referees
for helpful comments, and in particular for indicating that
Corollary~\ref{corkapparegular1} can be given in the present, more
general form.


%

\printaddresses

\end{document}